\documentclass[12pt,a4paper,reqno]{amsart}  
\usepackage{mathrsfs}
\usepackage{mathtools}
\usepackage{url}
 
\usepackage[headings]{fullpage}
\usepackage{comment} 
\usepackage{datetime}  
\usepackage[english]{babel} 

\usepackage[full]{textcomp}
\usepackage[osf]{newtxtext} 
\usepackage{cabin} 
\usepackage[bigdelims,vvarbb]{newtxmath} 
\usepackage[cal=boondoxo]{mathalfa} 
\usepackage{zlmtt}

\usepackage{microtype} 

\newcommand{\eps}{\varepsilon}

\newcommand{\dx}{\mathrm{d}}

\newcommand{\E}{\mathcal{E}}

\newcommand{\I}{\mathcal{I}} 
\newcommand{\J}{\mathcal{J}} 

\newcommand{\U}{\mathcal{U}}
\newcommand{\V}{\mathcal{V}}

\newcommand{\Scal}{\mathcal{S}}  

\newcommand{\C}{\mathbb{C}}
\newcommand{\N}{\mathbb{N}}
\newcommand{\R}{\mathbb{R}} 

\newcommand{\second}{{\prime\prime}}

\newcommand{\EE}{\mathscr{E}}
\newcommand{\Sgoth}{\mathscr{S}}

\newcommand{\Etilde}{\widetilde{E}}

\newcommand{\Stilde}{\widetilde{S}}
 
\newcommand{\Vtilde}{\widetilde{V}}

\newcommand{\Odi}[1]{\Odip{}{#1}}
\newcommand{\Odig}[1]{\mathcal{O}\Bigl(#1\Bigr)}
\newcommand{\Odim}[1]{\mathcal{O}\bigl(#1\bigr)}   
\newcommand{\Odip}[2]{\mathcal{O}_{#1}\left(#2\right)}
\newcommand{\Odipg}[2]{\mathcal{O}_{#1}\Bigl(#2\Bigr)}  
\newcommand{\Odipm}[2]{\mathcal{O}_{#1} (#2)}  
\newcommand{\odip}[2]{{o}_{#1}\left(#2\right)}
\newcommand{\odi}[1]{\odip{}{#1}}

\allowdisplaybreaks

\renewcommand{\qedsymbol}{$\square$}
\newenvironment{Proof}[1][Proof]{\par\noindent\textbf{#1.}~}
{\hfill\qedsymbol\smallskip\par}

\newtheoremstyle{sltheorems}
{10pt}
{6pt}
{\slshape}
{}
{\bfseries}
{.}
{.5em}
{\thmname{#1}\thmnumber{ #2}\thmnote{ (#3)}}

\theoremstyle{sltheorems} 
\newtheorem{Theorem}{Theorem}
\newtheorem{Lemma}{Lemma} 

\newtheoremstyle{remark}
{10pt}
{6pt}
{\rm} 
{}
{\bfseries}
{.}
{.5em}
{\thmname{#1}\thmnumber{ #2}\thmnote{ (#3)}}
 \theoremstyle{remark}

\makeatletter

\def\env@Biggcases{%
  \let\@ifnextchar\new@ifnextchar
  \Biggl\lbrace
  \def\arraystretch{1.2}%
  \array{@{}l@{\quad}l@{}}%
}
\makeatother

\allowdisplaybreaks

\begin{document}

\title[Asymptotic formulae for binary problems]{Short intervals asymptotic formulae for binary problems with prime powers, II} 


\author[]{Alessandro Languasco \lowercase{and} Alessandro Zaccagnini}

\date{\today, \currenttime}




\begin{abstract}
We improve some results in our paper \cite{LanguascoZ2018a} about the
asymptotic formulae in short intervals 
for the average number  of representations
of  integers of the forms $n=p_{1}^{\ell_1}+p_{2}^{\ell_2}$
and $n=p^{\ell_1} + m^{\ell_2}$, 
where $\ell_1, \ell_2\ge 2$ are fixed integers,
 $p,p_1,p_2$ are prime numbers and $m$ is an integer.
\end{abstract}

\maketitle

\bigskip
\section{Introduction}
Let $N$ be a sufficiently large integer and $1\le H \le N$.
In our recent papers   \cite{LanguascoZ2016b} and \cite{LanguascoZ2017c}
we provided   suitable asymptotic formulae in short intervals $[N, N + H]$ for the  
number of representation of an integer $n$ as a sum of a prime and 
a prime square, as a sum of a prime and a square, as the sum of two prime squares
or as a sum of a prime square and a square. 
To describe these results we need the following definitions.
Let $\ell_1, \ell_2\ge 1$ be integers, 
\begin{equation}
\label{density-def-c-def}
\lambda: =1/\ell_1+1/\ell_2 \quad \textrm{and}\quad
c(\ell_1,\ell_2):=
\frac{\Gamma(1/\ell_1)\Gamma(1/\ell_2)}{\ell_1\ell_2\Gamma(\lambda)}
=
c(\ell_2,\ell_1).
\end{equation}  
Using these notations we can say that our results in \cite{LanguascoZ2016b} and \cite{LanguascoZ2017c}
are about $\lambda=3/2$ and $\lambda=1$ while here we are interested in the case $\lambda< 1$.
We also recall that   Suzuki \cite{Suzuki2017b,Suzuki2017} has recently sharpened
our results in   \cite{LanguascoZ2017c} for the case $\lambda=3/2$.
In  \cite{LanguascoZ2018a} we were able to get non trivial  results on the case $\lambda<1$
but unfortunately in the unconditional case we were not able to address every possible combination of $\ell_1, \ell_2$.
The aim of this paper is to remove such limitations thus getting non trivial unconditional results for every $\ell_1, \ell_2\ge 2$ such that $\lambda<1$.
Moreover we will also improve a conditional result 
contained in \cite{LanguascoZ2018a} by extending
its uniformity range.
Such improvements follow from better estimates
of the error term involved in the main terms treatments
and from using a Tolev's lemma on a truncated mean-square average for the exponential sums over primes,
see Lemmas \ref{tolev-pulito-tilde} and \ref{tolev-pulito}. 
We recall here some definition already given in \cite{LanguascoZ2018a}. Let

\begin{equation}
\label{A-def} 
A=A(N,d) :=
 \exp \Bigl( d   \Bigl( \frac{\log N}{\log \log N} \Bigr)^{1/3} \Bigr),
\end{equation} 
where $d$ is a real parameter (positive or negative) chosen according
to need, and
\begin{equation}
\label{R-def}
R^{\second}_{\ell_1,\ell_2}(n)= 
\sum_ {p_{1}^{\ell_1}+p_{2}^{\ell_2}=n }
\log p_{1} \log p_{2}.
\end{equation}
The general shape of $A$ depends on the saving over the trivial bound
in the unconditional part of Lemma~\ref{LP-Lemma-gen}.
 In this case, due to the symmetry of the problem,
 we can assume that $2 \le \ell_1 \le \ell_2$.
We can now state
\begin{Theorem}
\label{thm-uncond}
 Let $N\ge 2$, $1\le H \le N$, 
 $2 \le \ell_1 \le \ell_2$ be integers and $\lambda<1$.
Then, for every $\eps>0$, there exists $C=C(\eps)>0$ such that
\begin{align*} 
\sum_{n=N+1}^{N+H} &
R^{\second}_{\ell_1,\ell_2}(n)
=
c(\ell_1,\ell_2)
H  N^{\lambda-1} 
+
\Odipm{\ell_1,\ell_2}{ H N^{\lambda-1} A(N, -C(\eps)) },
\end{align*}
uniformly for    $N^{ 1-5/(6\ell_2) +\eps}\le H \le N^{1-\eps}$,
where $\lambda$ and $c(\ell_1,\ell_2)$ are defined in \eqref{density-def-c-def}. 
\end{Theorem}

This should be compared with Theorem 1.1 of \cite{LanguascoZ2018a};
here the uniformity on $H$ is much larger so that Theorem \ref{thm-uncond}
is non-trivial for every choice of $2 \le \ell_1 \le \ell_2$ with $\lambda<1$.
We also remark that the uniformity level for  $H$ in Theorem \ref{thm-uncond} is  the  expected optimal one given the known density
estimates for the non trivial zeroes of the Riemann zeta function. 

Assuming the Riemann Hypothesis (RH)  holds 
we get a non-trivial result for    $\sum_{n = N+1}^{N + H} 
R^{\second}_{\ell_1,\ell_2}(n)$ uniformly for every $2\le \ell_1 \le \ell_2$
and $H$ in some range. 
We use throughout the paper the convenient notation $f=\infty(g)$
for $g=\odi{f}$.

\begin{Theorem}
\label{thm-RH}
 Let $N\ge 2$, $1\le H \le N$,  $2\le \ell_1\le \ell_2$ be integers, $\lambda<1$,
 and assume the Riemann Hypothesis holds.
Then
\begin{equation*} 
   \sum_{n = N+1}^{N + H} 
R^{\second}_{\ell_1,\ell_2}(n)
=
c(\ell_1,\ell_2)HN^{\lambda-1}
+\Odipm{\ell_1,\ell_2}{
H^2N^{\lambda-2}
+
H^{1/2} N^{1/\ell_1+1/(2\ell_2)-1/2}(\log N)^3
}
\end{equation*}
uniformly for  
 $\infty(N^{1-1/\ell_2}( \log N)^{6})\le H \le\odi{N}$,
 where $\lambda$ and $c(\ell_1,\ell_2)$ are defined in 
 \eqref{density-def-c-def}.
\end{Theorem}
This should be compared with Theorem 1.2 of \cite{LanguascoZ2018a}.
Here the second 
error term is improved and, as a consequence, the uniformity on $H$ is much larger
and essentially optimal given the spacing of the sequences.
If $\ell_1=2$ the log-power in the final result
can be slightly improved by using Lemma \ref{zac-lemma-series}  instead
of Lemma \ref{tolev-pulito-tilde} but, since the
improvement is marginal, we did not insert this
estimate in the proof of Theorem \ref{thm-RH}.

A slightly different problem is the one in which we replace
a prime power with a power. 
Letting
 \[
 r^{\prime}_{\ell_1,\ell_2}(n) = \!\!\!\!\!\!\!\!
\sum_{\substack {p^{\ell_1}+m^{\ell_2}=n\\N/A \le  p^{\ell_1} , \, m^{\ell_2} \le  N}}
\!\!\!\!\!\!\!\!
\log p,
\]
we lose the symmetry in $\ell_1 ,\ell_2$; hence we just assume 
that $\ell_1, \ell_2\ge 2$.
Here we take $A$ as defined in \eqref{A-def} with a suitable $d > 0$.
We need to change the setting and to use the 
finite sums, see Section~\ref{setting-finite-sums}, because, with the unique exception of the case $\ell=2$,
we cannot use the infinite series in this problem; this, for $\ell\ne 2$, is due to the 
lack of a suitable modular relation for the function
$ \omega_{\ell}(\alpha) =
\sum_{m=1}^{\infty}   
e^{-m^{\ell}/N} e(m^{\ell}\alpha) $, $\alpha\in [-1/2,1/2]$. For technical reasons, in this case we need to localise
the summands to get a sufficiently strong estimate in Lemma \ref{media-f-ell} below.

We have the following
\begin{Theorem}
\label{thm-uncond-HL}
 Let $N\ge 2$, $1\le H \le N$,  $\ell_1, \ell_2\ge 2$ be integers and $\lambda<1$.
Then, for every $\eps>0$, there exists $C=C(\eps)>0$ such that
\begin{align*} 
\sum_{n=N+1}^{N+H} &
r^{\prime}_{\ell_1,\ell_2}(n)
=
c(\ell_1,\ell_2)H  N^{\lambda-1} 
+
\Odipm{\ell_1,\ell_2}{H N^{\lambda-1} A(N, -C(\eps)) },
\end{align*}
uniformly for    $\max(N^{1-5/(6\ell_1)}; N^{1-1/\ell_2})N^{\eps}\le H \le N^{1-\eps}$,
where $\lambda$ and $c(\ell_1,\ell_2)$ are defined in \eqref{density-def-c-def}. 
\end{Theorem}
 
This should be compared with Theorem 1.3 of \cite{LanguascoZ2018a};
here the uniformity on $H$ is much larger so that Theorem \ref{thm-uncond-HL}
is non-trivial for every choice of $\ell_1,\ell_2\ge 2$ with $\lambda<1$.
We also remark that the uniformity level for  $H$ in Theorem \ref{thm-uncond-HL} is  the  expected optimal one given the known density
estimates for the non trivial zeroes of the Riemann zeta function
and the spacing of the sequences.

We finally remark that even assuming the Riemann Hypothesis we cannot improve 
the size of the error term in Theorem \ref{thm-uncond-HL} because in the main term evaluation we have
a term of the size $HN^{\lambda-1}  A^{- 1/\ell_2}$, see \eqref{main-term} below; 
moreover the magnitude of the error in the approximation in \eqref{T-ell-estim} 
is huge in the periphery of the arc, \emph{i.e.}, for $\alpha$ ``near'' $1/2$. So, under the assumption
of the Riemann Hypothesis, we can improve  Theorem \ref{thm-uncond-HL} essentially
 only for $\ell=2$ by using the infinite series approach; but this result is already
presented in \cite{LanguascoZ2018a}.

The basic strategy for all of the proofs of our results is the same.
We rewrite the quantity we are studying as a suitable integral of a
product of exponential sums.
We replace these by simpler approximations, and then evaluate the
``main term'' and estimate the error terms that arise in the
approximations by means of the Lemmas proved in the next Section.
The drawback of using finite sums instead of infinite series is that
the main term has a more complicated shape and its treatment is less
straightforward.
The main new ingredient, and the reason why we can improve our earlier
results, is a consequence of a result due to Tolev \cite{Tolev1992}:
we need the two variants for infinite series and finite sums which we
state as Lemmas \ref{tolev-pulito-tilde} and \ref{tolev-pulito}.
In the proofs of Theorems \ref{thm-uncond}-\ref{thm-RH} we also
exploit the stronger error we have in \eqref{I1-eval-series}, whereas
in the remaining proofs we use the $L^2$ bound provided by Lemma
\ref{media-f-ell} instead of the $L^{\infty}$ bound.

   \section{Setting and lemmas for Theorems \ref{thm-uncond}-\ref{thm-RH}}
   \label{setting-series}
Let $\ell,\ell_1, \ell_2\ge 2$ be  integers,  $e(\alpha) = e^{2\pi i\alpha}$,
$\alpha \in [-1/2,1/2]$.
For proving the first two theorems is convenient to use the original 
Hardy-Littlewood functions because the main term contribution can be easier
evaluated comparing with the setting with  finite exponential sums.
Let
\begin{equation} 
\label{tilde-main-defs}
\Stilde_\ell(\alpha) = \sum_{n=1}^{\infty} \Lambda(n) e^{-n^{\ell}/N} e(n^{\ell}\alpha),  \quad
\Vtilde_\ell(\alpha) = \sum_{p=2}^{\infty} \log p \, e^{-p^{\ell}/N} e(p^{\ell}\alpha),
\end{equation} 
\textrm{and}
\begin{equation*} 
z= 1/N-2\pi i\alpha.
\end{equation*} 

We now list some results we will use later.
The lemmas in this Section are mostly bounds for exponential sums
of various types.
We will use them in Section~\ref{sec:proof-th1}, after the dissection
of the unit interval into subintervals where different tools are
needed to evaluate the main term and estimate error terms.

 \begin{Lemma}[Lemma 3 of  \cite{LanguascoZ2016b}]
\label{tilde-trivial-lemma}
Let $\ell\ge 1$ be an integer. Then
\(
\vert \Stilde_{\ell}(\alpha)- \Vtilde_{\ell}(\alpha) \vert 
\ll_{\ell}
 N^{1/(2\ell)}  .
\)
\end{Lemma}

\begin{Lemma}[Lemma 2 of  \cite{LanguascoZ2016a}]
Let $\ell \ge 1$ be an integer, $N \ge 2$  and $\alpha\in [-1/2,1/2]$.
Then
\begin{equation*}
\Stilde_{\ell}(\alpha)  
= 
\frac{\Gamma(1/\ell)}{\ell z^{1/\ell}}
- 
\frac{1}{\ell}\sum_{\rho}z^{-\rho/\ell}\Gamma\Bigl(\frac{\rho}{\ell}\Bigr) 
+
\Odip{\ell}{1},
\end{equation*}
where $\rho=\beta+i\gamma$ runs over
the non-trivial zeros of $\zeta(s)$.
\end{Lemma}
\begin{Lemma} [Lemma 4 of \cite{LanguascoZ2016a}]
 \label{Laplace-formula}
Let $N$ be a positive integer and 
$\mu > 0$.
Then, uniformly for $n \ge 1$ and $X>0$, we have
\[
  \int_{-X}^{X} z^{-\mu} e(-n \alpha) \, \dx \alpha
  =
  e^{- n / N} \frac{n^{\mu - 1}}{\Gamma(\mu)}
  +
  \Odipg{\mu}{\frac{1}{nX^\mu}},
\]
 where $\Gamma$ is Euler's
function.
\end{Lemma}

\begin{Proof}
We remark that the proof is identical to the one of Lemma 4 of \cite{LanguascoZ2016a}
but in that case we just stated the lemma in the particular case $X=1/2$. Now we need its full strength
and hence, for completeness, we rewrite its proof. 
We start with the identity
\[
  \frac1{2 \pi}
  \int_{\R} \frac{e^{i D u}}{(a + i u)^s} \, \dx u
  =
  \dfrac{D^{s - 1} e^{- a D}}{\Gamma(s)},
\]
which is valid for $\sigma = \Re(s) > 0$ and $a \in \C$ with
$\Re(a) > 0$ and $D > 0$.
Letting $u = -2 \pi \alpha$ 
and taking $s = \mu$, $D = n$ and $a = N^{-1}$ we find
\[
  \int_{\R} \frac{e(- n \alpha)}{(N^{-1} - 2 \pi i \alpha)^\mu} \, \dx \alpha
  =
  \int_{\R} z^{- \mu} e(- n \alpha) \, \dx \alpha
  =
  \dfrac{n^{\mu - 1} e^{- n / N}}{\Gamma(\mu)}.
\]

For $0 < X < Y$ let
\[
  I(X, Y)
  =
  \int_X^Y \frac{e^{i D u}}{(a + i u)^\mu} \, \dx u.
\]
An integration by parts yields
\[
  I(X, Y)
  =
  \Bigl[
    \frac1{i D} \, \frac{e^{i D u}}{(a + i u)^\mu}
  \Bigr]_X^Y
  +
  \frac{\mu}D
  \int_X^Y \frac{e^{i D u}}{(a + i u)^{\mu + 1}} \, \dx u.
\]
Since $a > 0$, the first summand is $\ll_{\mu} D^{-1} X^{-\mu}$, uniformly.
The second summand is
\[
  \ll 
  \frac{\mu}D
  \int_X^Y \frac{\dx u}{u^{\mu + 1}}
  \ll_\mu
  D^{-1} X^{-\mu}.
\]
The result follows.
\end{Proof}
\begin{Lemma}[Lemma 3 of \cite{LanguascoZ2016a} and 
Lemma 1 of \cite{LanguascoZ2016b}]
 \label{LP-Lemma-gen} 
Let $\eps$ be an arbitrarily small
positive constant,  $\ell \ge 1$ be an integer, $N$ be a
sufficiently large integer and $L= \log N$. Then there exists a positive constant 
$c_1 = c_{1}(\eps)$, which does not depend on $\ell$, such that 
\[
\int_{-\xi}^{\xi} \,
\Bigl\vert
\Stilde_\ell(\alpha) - \frac{\Gamma(1/\ell)}{\ell z^{1/\ell}}
\Bigr\vert^{2}
\dx \alpha 
\ll_{\ell}
 N^{2/\ell-1} A(N, - c_{1}) 
\]
uniformly for $ 0\le \xi < N^{-1 +5/(6\ell) - \eps}$.
Assuming RH we get 
\[
\int_{-\xi}^{\xi} \,
\Bigl\vert
\Stilde_\ell(\alpha) - \frac{\Gamma(1/\ell)}{\ell z^{1/\ell}}
\Bigr\vert^{2}
\dx \alpha 
\ll_{\ell}
N^{1/\ell}\xi L^{2}
\]
uniformly  for  $0 \le \xi \le 1/2$.
\end{Lemma} 
Some of the following lemmas hold for a real index 
$k$ instead of an integral one $\ell$; in general we will always use
$k$ to denote a real index.
The new ingredient we are using here is based on a Tolev's lemma
\cite{Tolev1992}.
\begin{Lemma}[Tolev]
\label{tolev-pulito-tilde}
Let  $k>1$, $n\in \N$ and $\tau>0$. 
We have
\[
 \int_{-\tau}^{\tau}\vert \Stilde_k(\alpha)\vert ^2\, \dx \alpha
\ll_k
\bigl(\tau N^{1/k}+N^{2/k-1}\bigr) L^3
\quad
\textrm{and}
\quad
 \int_{-\tau}^{\tau}\vert \Vtilde_k(\alpha)\vert ^2\, \dx \alpha
\ll_k
\bigl(\tau N^{1/k}+N^{2/k-1}\bigr) L^3.
\]
\end{Lemma} 

\begin{Proof}
We  just prove the first part since the second one follows immediately by remarking that
 the   primes are supported on a thinner set than   the prime powers.
Let  $P=(2NL/k)^{1/k}$. A direct estimate
gives   $\Stilde_{k}(\alpha)= \sum_{n\le P} \Lambda(n) e^{-n^k/N} e(n^k\alpha) + \Odipm{k}{L^{1/k}}$.
Recalling that the Prime Number Theorem implies
$S_{k}(\alpha;t) := \sum_{n\le t} \Lambda(n)   e(n^{k}\alpha)  \ll t$, a
partial integration argument gives 
\[
\sum_{n\le P} \Lambda(n) e^{-n^k/N} e(n^k\alpha) 
= 
-\frac{k}{N} \int_1^P t^{k-1} e^{-t^k/N} S_{k}(\alpha;t) \  \dx t
+ \Odipm{k}{L^{1/k}}.
\]
Using the inequality $(\vert a\vert + \vert b \vert)^{2} \ll \vert a\vert^{2} + \vert b \vert^{2}$, Cauchy-Schwarz inequality
and interchanging the integrals, we get that
\begin{align*}
\int_{-\tau}^{\tau}
&\vert \Stilde_{k}(\alpha)\vert ^{2} \ \dx\alpha 
\ll_{k} 
\int_{-\tau}^{\tau}
\Bigl\vert \frac{1}{N} \int_1^P t^{k-1} e^{-t^k/N} S_{k}(\alpha;t)\  \dx t \Bigr\vert ^{2} \ \dx\alpha 
+  L^{2/k}
\\
&\ll_{k} 
\frac{1}{N^{2}}
\Bigl( \int_1^P t^{k-1} e^{-t^k/N} \  \dx t \Bigr)
\Bigl( \int_1^P t^{k-1} e^{-t^k/N} \int_{-\tau}^{\tau}\vert S_{k}(\alpha;t)\vert ^{2} \dx\alpha  \  \dx t \Bigr)
+  L^{2/k}.
\end{align*} 
Lemma 7 of Tolev \cite{Tolev1992} in the form given in Lemma 5 of \cite{GambiniLZ2018} on 
$S_{k}(\alpha;t) = \sum_{n\le t} \Lambda(n) e(n^k\alpha)$ implies
that 
$\int_{-\tau}^{\tau} \vert S_{k}(\alpha;t)\vert^{2}\ \dx\alpha \ll_{k}
\bigl(\tau t+t^{2-k}\bigr)(\log t)^3 $.
Using such an estimate 
and remarking that $ \int_1^P t^{k-1} e^{-t^k/N} \  \dx t \ll_k N$,
we obtain  that
\begin{align*}
\int_{-\tau}^{\tau}
\vert \Stilde_{k}(\alpha)\vert ^{2} \ \dx\alpha 
& \ll_{k} 
\frac{1}{N} 
 \int_1^P  \bigl(\tau t+t^{2-k}\bigr) t^{k-1}e^{-t^k/N} (\log t)^3 \  \dx t  
+  L^{2/k}
\\&
 \ll_{k}  
 \bigl(\tau N^{1/k}+N^{2/k-1}\bigr)L^{3} 
\end{align*}
by a direct computation.
This proves the first part of the lemma. 
\end{Proof}
In the case $\ell=2$ a slightly better final result
can be obtained using
\begin{Lemma}[Lemma 2 of \cite{LanguascoZ2017c}]
\label{zac-lemma-series}
Let $\ell\ge 2 $ be an integer and $0<\tau\leq 1/2$. Then
\[
\int_{-\tau}^{\tau} 
|\Stilde_{\ell}(\alpha)|^2 \ \dx\alpha 
\ll_{\ell} 
\tau N^{1/\ell} L  +
\begin{cases}
L^{2} & \text{if}\ \ell =2\\
1 & \text{if}\ \ell > 2.
\end{cases}
\]
\end{Lemma}

The next lemma is a consequence of Lemmas \ref{LP-Lemma-gen}-\ref{tolev-pulito-tilde}; its proof follows the line of Lemma \ref{tolev-coda}.
\begin{Lemma}
\label{tolev-coda-tilde}
Let $N\in \N$,   $k>1$,   $u\ge 1$ and $N^{-u}\le \omega \le N^{1/k-1}/L$. 
Let further $I(\omega):=[-1/2,-\omega]\cup  [\omega, 1/2]$.
We have
\[
\int_{I(\omega)}
\vert \Stilde_k(\alpha)\vert ^2 \frac{\dx \alpha}{\vert \alpha\vert}
\ll_{k}
\frac{ N^{2/k-1}}{\omega} L^3
 \quad
 \text{and}
 \quad
\int_{I(\omega)}
\vert \Vtilde_k(\alpha)\vert ^2 \frac{\dx \alpha}{\vert \alpha\vert}
\ll_{k}
\frac{ N^{2/k-1}}{\omega} L^3.
\]
Let further assume the Riemann Hypothesis, $\ell\ge 1$ be an integer and  $N^{-u}\le \eta \le 1/2$. Then
\[
\int_{I(\eta)}
\Bigl\vert
\Stilde_\ell(\alpha) - \frac{\Gamma(1/\ell)}{\ell z^{1/\ell}}
\Bigr\vert^{2}
 \frac{\dx \alpha}{\vert \alpha\vert}
\ll_{\ell}
  N^{1/\ell} L^3.
\]
\end{Lemma}
Let further
\[
U(\alpha,H) := \sum_{1\le  m\le  H}e(m\alpha).
\]
 We also have the usual numerically explicit inequality
\begin{equation}
\label{UH-estim}
\vert U(\alpha,H) \vert
\le
\min \bigl(H;  \vert \alpha\vert ^{-1}\bigr),
\end{equation}
see, \emph{e.g.}, on page 39 of Montgomery \cite{Montgomery1994}.
Using   \eqref{UH-estim} we obtain
 \begin{Lemma} 
%
Let $H\ge 2$, $\mu\in \R$, $\mu\ge 1$. Then
\begin{equation}
\label{U-media-def}
\U(\mu, H)
:=
\int_{-1/2}^{1/2} 
\vert U(\alpha,H)\vert^{\mu} \, \dx \alpha 
\ll
\begin{cases}
\log H& \text{if}\ \mu=1 \\
H^{\mu-1}& \text{if}\ \mu>1.
\end{cases}
\end{equation}
\end{Lemma}

Combining  \eqref{UH-estim},
Lemmas \ref{LP-Lemma-gen} and \ref{tolev-coda-tilde}
we get
\begin{Lemma}
Let   $\ell \ge 1$ be an integer, $N$ be a
sufficiently large integer and $L= \log N$.
Assume the Riemann Hypothesis. We have 
\begin{equation}
\label{E-media-def}
\E_{\ell}(H,N)
:=
\int_{-1/2}^{1/2} \,
\Bigl\vert
\Stilde_\ell(\alpha) - \frac{\Gamma(1/\ell)}{\ell z^{1/\ell}}
\Bigr\vert^{2}
\vert U(-\alpha,H) \vert
\ \dx \alpha 
\ll_{\ell}
N^{1/\ell}  L^{3}.
\end{equation}
\end{Lemma}

Combining  \eqref{UH-estim},
Lemmas \ref{tolev-pulito-tilde}, \ref{tolev-coda-tilde} and \ref{tilde-trivial-lemma}
we get
\begin{Lemma}
Let   $k > 1$, $N$ be a
sufficiently large integer, $L= \log N$ and    $N^{1-1/k}L \ll H \le N$.  We have 
\begin{equation}
\label{S-media-def}
\Scal_{k}(H,N)
:=
\int_{-1/2}^{1/2}  
\vert
\Stilde_k(\alpha)
\vert^{2}
\vert U(-\alpha,H) \vert
\ \dx \alpha 
\ll_{k}
H N^{2/k-1}  L^{3}
\end{equation}
and
\begin{equation}
\label{V-media-def}
\V_{k}(H,N)
:=
\int_{-1/2}^{1/2}  
\vert
\Vtilde_k(\alpha)
\vert^{2}
\vert U(-\alpha,H) \vert
\ \dx \alpha 
\ll_{k}
H N^{2/k-1}  L^{3}.
\end{equation}
\end{Lemma}
   
 \section{Proof of Theorem \ref{thm-uncond}}
\label{sec:proof-th1}
Due to the symmetry of the summands we may let $2\le \ell_1\le \ell_2$, $\lambda<1$,
where  $\lambda$ is defined in 
 \eqref{density-def-c-def}; 
we'll see at the end
of the proof how the conditions in the statement of this theorem
follow.
Assume 
\begin{equation}
\label{B-def-th1}
B=B(N,\eps)=
N^{\eps}
\end{equation} 
and let $H>2B$.
Basically, we now replace $\Vtilde_{\ell}$ by $\Stilde_{\ell}$ at the
centre of the integration interval, that is on $[-B/H, B/H]$.
Then we bound the error term and the contribution of the remainder of
the integration range by means of several Lemmas proved in
Section~\ref{setting-series}.
We have 
\begin{align}
\notag 
\sum_{n=N+1}^{N+H} e^{-n/N} 
R^{\second}_{\ell_1,\ell_2}(n) 
&= 
\int_{-1/2}^{1/2} \Vtilde_{\ell_1}(\alpha)  \Vtilde_{\ell_2}(\alpha)  U(-\alpha,H)e(-N\alpha) \, \dx \alpha
\\
\notag
 &=
\int_{-B/H}^{B/H} \Stilde_{\ell_1}(\alpha) \Stilde_{\ell_2}(\alpha)  U(-\alpha,H)e(-N\alpha) \, \dx \alpha 
\\
\notag
&
\hskip1cm+
\int_{I(B/H)} 
\Stilde_{\ell_1}(\alpha) \Stilde_{\ell_2}(\alpha)  U(-\alpha,H)e(-N\alpha) \, \dx \alpha
\\
\notag
&
\hskip1cm+
\int_{-1/2}^{1/2} (\Vtilde_{\ell_1}(\alpha)-\Stilde_{\ell_1}(\alpha)) \Vtilde_{\ell_2}(\alpha)  U(-\alpha,H)e(-N\alpha) \, \dx \alpha 
\\
\notag
&
\hskip1cm+
\int_{-1/2}^{1/2} (\Vtilde_{\ell_2}(\alpha)-\Stilde_{\ell_2}(\alpha)) \Stilde_{\ell_1}(\alpha)  U(-\alpha,H)e(-N\alpha) \, \dx \alpha,
\\&
\label{dissect} 
=I_1+I_2+I_3+I_4,
\end{align}
say,
where  $I(B/H):=[-1/2,-B/H]\cup  [B/H, 1/2]$.

\subsection{Estimate of $I_2$}
By the Cauchy-Schwarz inequality, \eqref{UH-estim}  and Lemma \ref{tolev-coda-tilde} we have
\begin{equation} 
  \label{S-average-coda} 
I_2 
 \ll
\Bigl(
\int_{I(B/H)} 
\vert \Stilde_{\ell_1}(\alpha) \vert ^{2}  \frac{\dx \alpha}{\vert \alpha \vert}  
\Bigr)^{1/2} 
\Bigl(
\int_{I(B/H)} 
\vert \Stilde_{\ell_2}(\alpha) \vert ^{2} \frac{\dx \alpha}{\vert \alpha \vert} 
\Bigr)^{1/2}
\ll_{\ell_1,\ell_2}
  \frac{H N^{\lambda -1}L^3}{B},
\end{equation}
provided that $H\gg \max( N^{1-1/\ell_1}; N^{1-1/\ell_2})BL=N^{1-1/\ell_2}BL$,
since $\ell_2\ge \ell_1$.

\subsection{Estimate of $I_3$ and $I_4$}
\label{sec:I3-I4-estim}
%
Using Lemma \ref{tilde-trivial-lemma}, 
the Cauchy-Schwarz inequality,
\eqref{V-media-def} and
\eqref{U-media-def},
we get
\begin{equation}  
\label{I3-estim}
I_3
 \ll_{\ell_1,\ell_2}
N^{1/(2\ell_1)}
\V_{\ell_2}(H,N)^{1/2}
\U(1,H)^{1/2}
\ll_{\ell_1,\ell_2}
H^{1/2}N^{1/(2\ell_1)+1/\ell_2-1/2}L^2,
\end{equation}
provided that $H\gg N^{1-1/\ell_2}L$.
Using Lemma \ref{tilde-trivial-lemma}, 
the Cauchy-Schwarz inequality,
\eqref{S-media-def} and
\eqref{U-media-def},
we get
\begin{equation}  
\label{I4-estim}
I_4
\ll_{\ell_1,\ell_2}
N^{1/(2\ell_2)}
\Scal_{\ell_1}(H,N)^{1/2}
\U(1,H)^{1/2}  
 \ll_{\ell_1,\ell_2}
H^{1/2}N^{1/\ell_1+1/(2\ell_2)-1/2}L^2.
\end{equation}
provided that $H\gg N^{1-1/\ell_1}L$. Hence, using \eqref{I3-estim}-\eqref{I4-estim}
and recalling $\ell_1\le \ell_2$, we have
\begin{equation}  
\label{I3-I4-estim}
I_3+I_4 
\ll_{\ell_1,\ell_2}
H^{1/2}N^{1/\ell_1+1/(2\ell_2)-1/2}L^2.
\end{equation}
provided that $H\gg N^{1-1/\ell_2}L$.

\subsection{Evaluation of  $I_1$}
We now obtain the main term.
From now on, we denote 
%
\[
\Etilde_{\ell}(\alpha) : =\Stilde_\ell(\alpha) - \frac{\Gamma(1/\ell)}{\ell z^{1/\ell}}.
\]
%
In the formula below, we see that the main term arises from
the product of the two terms $\Gamma(1 / \ell) / (\ell z^{1/\ell})$.
The other terms give a smaller contribution, since they contain at
least one factor $\Etilde_{\ell}$, which is small on average by
Lemma~\ref{LP-Lemma-gen}.
Recallining \eqref{dissect},
by \eqref{tilde-main-defs} we get
\begin{align} 
\notag 
I_1  
&=
\frac{\Gamma(1/\ell_1)\Gamma(1/\ell_2)}{\ell_1\ell_2}
\int_{-B/H}^{B/H} z^{-1/\ell_1-1/\ell_2}  U(-\alpha,H)e(-N\alpha) \, \dx \alpha 
\\
\notag
&
\hskip0.5cm
+
\frac{\Gamma(1/\ell_1)}{\ell_1}
\int_{-B/H}^{B/H} z^{-1/\ell_1} \Etilde_{\ell_2}(\alpha)
U(-\alpha,H)e(-N\alpha) \, \dx \alpha 
\\
 \notag
 &
\hskip0.5cm
+
\frac{\Gamma(1/\ell_2)}{\ell_2}
\int_{-B/H}^{B/H} z^{-1/\ell_2}\Etilde_{\ell_1}(\alpha) U(-\alpha,H)e(-N\alpha) \, \dx \alpha 
\\
 \notag
 &
\hskip0.5cm
+
\int_{-B/H}^{B/H} \Etilde_{\ell_1}(\alpha)\Etilde_{\ell_2}(\alpha)
 U(-\alpha,H)e(-N\alpha) \, \dx \alpha  
\\
 \label{main-dissection}
&= \I_1 +\I_2 + \I_3 + \I_4  ,
 \end{align}
 say.
 We now evaluate these terms.
 
\subsection{Computation of the main term $\I_1$}
 
By Lemma \ref{Laplace-formula},
 \eqref{density-def-c-def}
and using  $e^{-n/N}=e^{-1}+ \Odi{H/N}$ for $n\in[N+1,N+H]$, $1\le H \le N$,
a direct calculation  gives
\begin{align}
\notag
\I_1
&
=
c(\ell_1,\ell_2)
\sum_{n=N+1}^{N+H} e^{-n/N} n^{\lambda-1}
+\Odipg{\ell_1,\ell_2}{\frac{H}{N}\Bigl(\frac{H}{B}\Bigr)^{\lambda}}
\\
\notag
&
=
\frac{c(\ell_1,\ell_2)}{e}
\sum_{n=N+1}^{N+H} n^{\lambda-1}
+\Odipg{\ell_1,\ell_2}{\frac{H}{N}\Bigl(\frac{H}{B}\Bigr)^{\lambda}+H^2N^{\lambda-2}}
\\
\label{I1-eval-series}
&
=
c(\ell_1,\ell_2)
\frac{HN^{\lambda-1}}{e}
+\Odipg{\ell_1,\ell_2}{\frac{H}{N}\Bigl(\frac{H}{B}\Bigr)^{\lambda}+H^2N^{\lambda-2}+N^{\lambda-1}}.
 \end{align}

\subsection{Estimate of $\I_4$}
Denote
\begin{equation}
\label{E-media-def-unc}
\EE_{\ell}(B/H,N)
:=
\int_{-B/H}^{B/H}
\vert  \Etilde_{\ell}(\alpha) \vert^2 \, \dx \alpha
\ll_\ell
N^{2/\ell-1}A(N, -c_1),
\end{equation}
in which the estimate follows from
Lemma \ref{LP-Lemma-gen} 
provided that $H\gg N^{1-5/(6\ell)+\eps}B$.
 Using \eqref{density-def-c-def},
  the Cauchy-Schwarz inequality
 and \eqref{E-media-def-unc},
 we obtain
 \begin{align}
\notag
\I_4
& \ll_{\ell_1,\ell_2}
H
\int_{-B/H}^{B/H}
\vert \Etilde_{\ell_1}(\alpha) \vert \vert \Etilde_{\ell_2}(\alpha) \vert \, \dx \alpha
\ll_{\ell_1,\ell_2}
H
\EE_{\ell_1}(B/H,N)^{1/2}
\EE_{\ell_2}(B/H,N)^{1/2}
\\
\label{I4-estim-series}&
\ll_{\ell_1,\ell_2} 
H  N^{\lambda-1} A(N, -c_1/2),
\end{align}
provided that $H\gg N^{1-5/(6\ell_2)+\eps}B$.

\subsection{Estimate of $\I_2$}
Denote
\begin{equation}
\label{S-media-def-unc}
\Sgoth_{\ell}(B/H,N)
:=
\int_{-B/H}^{B/H}
\vert  \Stilde_{\ell}(\alpha) \vert^2 \, \dx \alpha
\ll_\ell
N^{2/\ell-1} L^3,
\end{equation}
in which the estimates follow from
Lemma   \ref{tolev-pulito},
provided that $H\gg N^{1-1/\ell}B$.
Remarking 
$\vert z\vert^{-1/\ell} 
\ll_{\ell}
\vert \Stilde_{\ell}(\alpha) \vert
+
\vert \Etilde_{\ell}(\alpha) \vert,
$
using the Cauchy-Schwarz inequality, 
  \eqref{E-media-def-unc}
  and \eqref{S-media-def-unc},
we obtain
\begin{align}
\notag
\I_2 &\ll_{\ell_1,\ell_2}
H
\int_{-B/H}^{B/H} 
\vert \Stilde_{\ell_1}(\alpha) \vert
\vert \Etilde_{\ell_2}(\alpha) \vert
 \, \dx \alpha 
+
H
\int_{-B/H}^{B/H} 
\vert \Etilde_{\ell_1}(\alpha) \vert
\vert \Etilde_{\ell_2}(\alpha) \vert
\, \dx \alpha 
\\
\notag
&\ll_{\ell_1,\ell_2}
H
\Sgoth_{\ell_1}(B/H,N)^{1/2}
\EE_{\ell_2}(B/H,N)^{1/2}
+ 
H
\EE_{\ell_1}(B/H,N)^{1/2}
\EE_{\ell_2}(B/H,N)^{1/2}
\\
\label{I2-estim-unc}
&
\ll_{\ell_1,\ell_2}
H  N^{\lambda-1} A(N, -c_1/4),
\end{align}
provided that  $N^{1-5/(6\ell_2)+\eps} B \le H \le N^{1-\eps}$.

\subsection{Estimate of $\I_3$}
It's very similar to $\I_2$'s;
we just need to interchange $\ell_1$ with $\ell_2$
thus getting that there exists $C=C(\eps)>0$ such that 
 \begin{equation}
\label{I3-estim-series}
\I_3 \ll_{\ell_1,\ell_2}
 H N^{\lambda-1} A(N, -C),
\end{equation}
provided that  $N^{1-5/(6\ell_2)+\eps} B \le H \le N^{1-\eps}$.

\subsection{Final words}
 
Summarizing, recalling that $2\le\ell_1\le \ell_2$, $\lambda<1$, by
\eqref{dissect}-\eqref{I4-estim}, 
 \eqref{main-dissection}-\eqref{I1-eval-series}, 
\eqref{I4-estim-series},
\eqref{I2-estim-unc}-\eqref{I3-estim-series} and
by  optimising
the choice of $B$ as in \eqref{B-def-th1},
we have that there exists $C=C(\eps)>0$ such that 
\begin{align}
\label{almost-done-th1} 
\sum_{n=N+1}^{N+H} e^{-n/N} 
R^{\second}_{\ell_1,\ell_2}(n) 
=
\frac{c(\ell_1,\ell_2)}{e}   H N^{\lambda-1} 
+
\Odipm{\ell_1,\ell_2}{H N^{\lambda-1} A(N, -C)},
\end{align}
uniformly for $N^{1-5/(6\ell_2) +\eps}\le H \le N^{1-\eps}$.
{}From  $e^{-n/N}=e^{-1}+ \Odi{H/N}$ for $n\in[N+1,N+H]$, $1\le H \le N$,
we get 
\begin{align*}
   \sum_{n = N+1}^{N + H} 
R^{\second}_{\ell_1,\ell_2}(n)
&=
c(\ell_1,\ell_2)HN^{\lambda-1}
+
\Odipm{\ell_1,\ell_2}{H N^{\lambda-1} A(N, -C)}
+
  \Odig{\frac{H}{N}\sum_{n = N+1}^{N + H} R^{\second}_{\ell_1,\ell_2}(n)
}.
\end{align*}
Using $e^{n/N}\le  e^{2}$ 
and \eqref{almost-done-th1},
 the last error term is
$\ll_{\ell_1,\ell_2} H^2N^{\lambda-2}$.
Hence we get
\begin{equation*}
   \sum_{n = N+1}^{N + H} 
R^{\second}_{\ell_1,\ell_2}(n)
=
c(\ell_1,\ell_2)HN^{\lambda-1}
+\Odipm{\ell_1,\ell_2}{
 H N^{\lambda-1} A(N, -C)
},
\end{equation*}
uniformly for $N^{1-5/(6\ell_2) +\eps}\le H \le N^{1-\eps}$ and $2\le \ell_1\le \ell_2$.
Theorem \ref{thm-uncond} follows.

\section{Proof of Theorem \ref{thm-RH}}
In this Section we assume the Riemann Hypothesis holds.
In this case, we do not need a different argument for ``centre'' and
``periphery'' of the integration interval, since
Lemma~\ref{LP-Lemma-gen} is valid throughout $[-1/2, 1/2]$.
Recalling \eqref{R-def},
we have  
\begin{align}
\notag
\sum_{n=N+1}^{N+H} e^{-n/N} 
R^{\second}_{\ell_1,\ell_2}(n)
&= 
\int_{-1/2}^{1/2} \Vtilde_{\ell_1}(\alpha) \Vtilde_{\ell_2}(\alpha)  U(-\alpha,H)e(-N\alpha) \, \dx \alpha 
\\
\notag
&
=
\int_{-1/2}^{1/2} \Stilde_{\ell_1}(\alpha) \Stilde_{\ell_2}(\alpha)  U(-\alpha,H)e(-N\alpha) \, \dx \alpha 
\\ 
\notag
&
\hskip1cm+
\int_{-1/2}^{1/2} (\Vtilde_{\ell_1}(\alpha)-\Stilde_{\ell_1}(\alpha)) \Vtilde_{\ell_2}(\alpha)  U(-\alpha,H)e(-N\alpha) \, \dx \alpha 
\\
\notag
&
\hskip1cm+
\int_{-1/2}^{1/2} (\Vtilde_{\ell_2}(\alpha)-\Stilde_{\ell_2}(\alpha)) \Stilde_{\ell_1}(\alpha)  U(-\alpha,H)e(-N\alpha) \, \dx \alpha,
\\&
\label{dissect-RH} 
=J_1+J_2+J_3,
\end{align}
say.
\subsection{Estimate of $J_2$ and $J_3$}
The quantities $J_2$ and $J_3$ are equal to $I_3$ and $I_4$ of Section~\ref{sec:I3-I4-estim}.
Hence by \eqref{I3-I4-estim} we get
\begin{equation}  
\label{J2-J3-estim}
J_2+J_3 
\ll_{\ell_1,\ell_2}
H^{1/2}N^{1/\ell_1+1/(2\ell_2)-1/2}L^2.
\end{equation}
provided that $H\gg N^{1-1/\ell_2}L$.

\subsection{Evaluation of  $J_1$}
Here we obtain the main term essentially as above, but we can deal
with the whole integration interval at once.
Hence
\begin{align}
 \notag
J_1
& = \frac{\Gamma(1/\ell_1)\Gamma(1/\ell_2)}{\ell_1\ell_2}
 \int_{-1/2}^{1/2} z^{-1/\ell_1-1/\ell_2}  U(-\alpha,H)e(-N\alpha) \, \dx \alpha 
\\
 \notag
 &
\hskip0.5cm
+
\frac{\Gamma(1/\ell_1)}{\ell_1}
\int_{-1/2}^{1/2} z^{-1/\ell_1} \Etilde_{\ell_2}(\alpha) 
U(-\alpha,H)e(-N\alpha) \, \dx \alpha 
\\
 \notag
 &
\hskip0.5cm
+
\frac{\Gamma(1/\ell_2)}{\ell_2}
\int_{-1/2}^{1/2}  z^{-1/\ell_2} \Etilde_{\ell_1}(\alpha)   U(-\alpha,H)e(-N\alpha) \, \dx \alpha 
\\
 \notag
 &
\hskip0.5cm
+
\int_{-1/2}^{1/2}  \Etilde_{\ell_1}(\alpha) \Etilde_{\ell_2}(\alpha) 
 U(-\alpha,H)e(-N\alpha) \, \dx \alpha 
\\
 \label{main-dissection-RH-series}
 &  
= \J_1 +\J_2 + \J_3 + \J_4,
 \end{align}
 say.
 Now we evaluate these terms.
 
\subsection{Computation of $\J_1$}

By Lemma \ref{Laplace-formula},
 \eqref{density-def-c-def}
and using  $e^{-n/N}=e^{-1}+ \Odi{H/N}$ for $n\in[N+1,N+H]$, $1\le H \le N$,
a direct calculation  gives
\begin{align}
\notag
\J_1
&
=
c(\ell_1,\ell_2)
\sum_{n=N+1}^{N+H} e^{-n/N} n^{\lambda-1}
+\Odipg{\ell_1,\ell_2}{\frac{H}{N}}
\\
\notag
&
=
\frac{c(\ell_1,\ell_2)}{e}
\sum_{n=N+1}^{N+H} n^{\lambda-1}
+\Odipg{\ell_1,\ell_2}{\frac{H}{N}+H^2N^{\lambda-2}}
\\
\label{J1-eval-RH-series}
&
=
c(\ell_1,\ell_2)
\frac{HN^{\lambda-1}}{e}
+\Odipg{\ell_1,\ell_2}{\frac{H}{N}+H^2N^{\lambda-2}+N^{\lambda-1}}.
 \end{align}

\subsection{Estimate of $\J_4$}
 Hence, by the Cauchy-Schwarz inequality and \eqref{E-media-def},
we obtain
 \begin{equation}
\label{J4-estim-RH-series}
\J_4
\ll_{\ell_1,\ell_2}
\E_{\ell_1}(H,N)^{1/2}
\E_{\ell_2}(H,N)^{1/2}
\ll_{\ell_1,\ell_2} 
  N^{\lambda/2} L^3.
\end{equation}

\subsection{Estimate of $\J_2$} %
Remarking 
$\vert z\vert^{-1/\ell} 
\ll_{\ell}
\vert \Stilde_{\ell}(\alpha) \vert
+
\vert \Etilde_{\ell}(\alpha) \vert,
$
using the Cauchy-Schwarz inequality, 
  \eqref{E-media-def}
  and \eqref{S-media-def},
we obtain

\begin{align}
\notag
\J_2 &\ll_{\ell_1,\ell_2}
\int_{-1/2}^{1/2} 
\vert \Stilde_{\ell_1}(\alpha) \vert
\vert \Etilde_{\ell_2}(\alpha) \vert
\vert 
U(-\alpha,H)
\vert \, \dx \alpha 
+
\int_{-1/2}^{1/2} 
\vert \Etilde_{\ell_1}(\alpha) \vert
\vert \Etilde_{\ell_2}(\alpha) \vert
\vert 
U(-\alpha,H)
\vert \, \dx \alpha 
\\
\notag
&\ll_{\ell_1,\ell_2}
\Scal_{\ell_1}(H,N)^{1/2}
\E_{\ell_2}(H,N)^{1/2}
+ 
\E_{\ell_1}(H,N)^{1/2}
\E_{\ell_2}(H,N)^{1/2}
\\
 \label{J2-estim-RH-series}  
 &
\ll_{\ell_1,\ell_2}
H^{1/2} N^{1/\ell_1+1/(2\ell_2)-1/2}L^3 ,
\end{align}
provided that $H\gg N^{1-1/\ell_1}L$.
  
\subsection{Estimate of $\J_3$}

The estimate of $\J_3$
is very similar to $\J_2$'s;
we just need to interchange $\ell_1$ with $\ell_2$.
We obtain
 \begin{equation}
\label{J3-estim-RH-series}
\J_3 \ll_{\ell_1,\ell_2}
H^{1/2} N^{1/\ell_2+1/(2\ell_1)-1/2}L^3,
\end{equation}
provided that $H\gg N^{1-1/\ell_2}L$.

\subsection{Final words}

Summarizing, recalling $2\le\ell_1\le \ell_2$, 
by  \eqref{density-def-c-def} and \eqref{dissect-RH}-\eqref{J3-estim-RH-series}
we have
\begin{align}
\notag
\sum_{n=N+1}^{N+H} 
e^{-n/N}
R^{\second}_{\ell_1,\ell_2}(n)
&=
c(\ell_1,\ell_2)
\frac{HN^{\lambda-1}}{e}
\\&
\label{almost-done}
+\Odipg{\ell_1,\ell_2}{
\frac{H}{N}
+
H^2N^{\lambda-2}
+
H^{1/2} N^{1/\ell_1+1/(2\ell_2)-1/2}L^3
}
\end{align}
which is an asymptotic formula 
for $\infty(N^{1-1/\ell_2}L^{6})\le H \le \odi{N}$.
{}From  $e^{-n/N}=e^{-1}+ \Odi{H/N}$ for $n\in[N+1,N+H]$, $1\le H \le N$,
we get 
\begin{align*}
   \sum_{n = N+1}^{N + H} 
R^{\second}_{\ell_1,\ell_2}(n)
&=
c(\ell_1,\ell_2)HN^{\lambda-1}
+\Odipg{\ell_1,\ell_2}{
H^2N^{\lambda-2}
+
H^{1/2} N^{1/\ell_1+1/(2\ell_2)-1/2}L^3
}  
\\&
+
  \Odig{\frac{H}{N}\sum_{n = N+1}^{N + H} R^{\second}_{\ell_1,\ell_2}(n)
}.
\end{align*}
Using $e^{n/N}\le  e^{2}$ 
and \eqref{almost-done},
 the last error term is
$\ll_{\ell_1,\ell_2} H^2N^{\lambda-2}$.
Hence we get
\begin{equation*}
   \sum_{n = N+1}^{N + H} 
R^{\second}_{\ell_1,\ell_2}(n)
=
c(\ell_1,\ell_2)HN^{\lambda-1}
+\Odipg{\ell_1,\ell_2}{
H^2N^{\lambda-2}
+
H^{1/2} N^{1/\ell_1+1/(2\ell_2)-1/2}L^3
},
\end{equation*}
uniformly for every $2\le \ell_1\le \ell_2$ and 
 $\infty(N^{1-1/\ell_2}L^{6})\le H \le\odi{N}$.
Theorem \ref{thm-RH} follows.

   \section{Setting and lemmas for Theorem \ref{thm-uncond-HL}}
   \label{setting-finite-sums}
 We also need similar lemmas for the finite sums since 
 we will use them for proving the second two results.
Let $k >0$ be a real number and
\begin{align}
\notag
S_{k}(\alpha) &:= \sum_{N/A \le  m^{k} \le  N}\Lambda(m)\ e(m^{k} \alpha),   \quad
V_{k}(\alpha) := \sum_{N/A \le  p^{k} \le  N} \log p\ e(p^{k} \alpha), \\
\label{main-defs}
T_{k}(\alpha) &:= \sum_{N/A \le m^{k} \le  N} e(m^{k} \alpha ), \quad
f_{k}(\alpha)  :=1/k\sum_{N/A \le  m\le  N} m^{1/k-1}\ e(m\alpha), 
\end{align}
where $A$ is defined in \eqref{A-def}.
As we remarked earlier, we take $d > 0$ in the definition of $A$.
We need this parameter because, if we chose $A = N$ in the definition
of $f_k$ above, the $L^2$ bound in Lemma~\ref{media-f-ell} would
become too weak.
We remark that we can choose $d$ in such a way that the constant
$C(\eps)$ in the statement of Theorem \ref{thm-uncond-HL} is
independent of $\ell_1$ and $\ell_2$.
By Lemmas 2.8 and 4.1 of Vaughan \cite{Vaughan1997}, we obtain
\begin{equation}
\label{T-ell-estim}
\vert T_{k}(\alpha) - f_{k}(\alpha) \vert \ll_k  (1+\vert \alpha \vert N)^{1/2}.
\end{equation}
We recall  that $\eps>0$ and we let $L= \log N$.
Now we recall some lemmas from \cite{LanguascoZ2018a}.

\begin{Lemma}[Lemma 2 of \cite{LanguascoZ2018a}]
%
Let $k >0$ be a real number. Then
\(
\vert S_{k}(\alpha)- V_{k}(\alpha) \vert 
\ll_{\ell}
 N^{1/(2k)}   .
\)
\end{Lemma} 

We need the following lemma which collects the results 
of Theorems 3.1-3.2 of \cite{LanguascoZ2013b}; see also Lemma 1 of \cite{LanguascoZ2016a}.
\begin{Lemma}   
Let $k  > 0$ be a real number and $\eps$ be an arbitrarily small
positive constant. Then there exists a positive constant 
$c_1 = c_{1}(\eps)$, which does not depend on $k$, such that
\[
\int_{-1/K}^{1/K}
\vert
S_{k}(\alpha) - T_{k}(\alpha)  
\vert^2 
\, \dx \alpha
\ll_{k}   N^{2/k-1}
\Bigl(
A(N, - c_{1})
+
\frac{K L^{2}}{N}
\Bigr),
\]
uniformly for  $N^{1-5/(6k)+\eps}\le K \le N$.  
Assuming further RH we get 
\[
\int_{-1/K}^{1/K}
\vert
S_{k}(\alpha) - T_{k}(\alpha) 
\vert^2 
\, \dx \alpha
\ll_{k} 
\frac{N^{1/k} L^{2}}{K} + K N^{2/k-2} L^{2},
\]
uniformly for  $N^{1-1/k}\le K \le N$.  
\end{Lemma} 
 
 Combining the two previous lemmas we get
\begin{Lemma}[Lemma 4 of \cite{LanguascoZ2018a}]
\label{App-BCP-Gallagher-2}
Let $k  > 0$ be a real number and $\eps$ be an arbitrarily small
positive constant. Then there exists a positive constant 
$c_1 = c_{1}(\eps)$, which does not depend on $k$, such that
\[
E_{k}(N,K) :=
\int_{-1/K}^{1/K}
\vert
V_{k}(\alpha) - T_k(\alpha)  
\vert^2 
\, \dx \alpha
\ll_{\ell}   N^{2/k-1}
\Bigl(
A(N, - c_{1})
+
\frac{K L^{2}}{N}
\Bigr),
\]
uniformly for  $N^{1-5/(6k)+\eps}\le K \le N$.  
Assuming further RH we get 
\[
E_{k}(N,K) 
\ll_{k} 
\frac{N^{1/k} L^{2}}{K} + K N^{2/k-2} L^{2},
\]
uniformly for  $N^{1-1/k}\le K \le N$.  
\end{Lemma}   

 \begin{Lemma}[Lemma 6 of \cite{LanguascoZ2018a}]
 \label{media-f-ell}
Let $k >0$ be a real number and recall that $A$ is defined in \eqref{A-def}. Then
\[
 F_{k}(N,A):=
\int_{-1/2}^{1/2} \vert f_{k}(\alpha) \vert ^2\, \dx \alpha 
\ll_{k} 
N^{2/k-1}
\begin{cases}
A^{1-2/k} & \text{if}\ k > 2\\
\log A & \text{if}\ k = 2\\
1 & \text{if}\ 0<k < 2.
\end{cases}
\]
 \end{Lemma} 

The new ingredient we are using here is based on a Tolev's lemma
\cite{Tolev1992}  in the form given in Lemma 5 of \cite{GambiniLZ2018}. 
\begin{Lemma}[Tolev]
\label{tolev-pulito}
Let  $k>1$, $n\in \N$ and $\tau>0$. Then
\[
 \int_{-\tau}^{\tau}\vert V_k(\alpha)\vert ^2\, \dx \alpha
\ll_k
\bigl(\tau N^{1/k}+N^{2/k-1}\bigr) L^3
\quad
\textrm{and}
\quad
\int_{-\tau}^{\tau}\vert T_k(\alpha)\vert ^2\, \dx \alpha
\ll_k
\bigl(\tau N^{1/k}+N^{2/k-1}\bigr) L.
\]
\end{Lemma} 

The last lemma is a consequence of Lemma \ref{tolev-pulito}.
\begin{Lemma}
\label{tolev-coda}
Let $N\in \N$,   $k>1$,   $c\ge 1$ and $N^{-c}\le \omega \le N^{1/k-1}/L$. Let further $I(\omega):=
[-1/2,-\omega]\cup[\omega,1/2]$.
Then we have
\[
\int_{I(\omega)}
\vert V_{k}(\alpha)\vert ^2 \frac{\dx \alpha}{\vert \alpha\vert}
\ll_{k}
\frac{ N^{2/k-1}}{\omega} L^3
 \quad
 \text{and}
 \quad
\int_{I(\omega)}
\vert T_{k}(\alpha)\vert ^2 \frac{\dx \alpha}{\vert \alpha\vert}
\ll_{k}
\frac{ N^{2/k-1}}{\omega} L.
\]
\end{Lemma}
\begin{Proof}
By partial integration and Lemma \ref{tolev-pulito} we get that
\begin{align*}
\int_{\omega}^{1/2}
\vert V_{k}(\alpha)\vert ^2 \frac{\dx \alpha}{ \alpha}
&\ll
\frac{1}{\omega}
\int_{-\omega}^{\omega} 
\vert V_{k}(\alpha) \vert^{2}\ \dx \alpha 
+
\int_{-1/2}^{1/2} 
\vert V_{k}(\alpha) \vert^{2}\ \dx \alpha 
+
\int_{\omega}^{1/2} 
\Bigl(
\int_{-\xi}^{\xi} 
\vert V_{k}(\alpha) \vert^{2}\ \dx \alpha 
\Bigr)
\frac{\dx \xi}{\xi^2}
\\&
\ll_{k}
\frac{1}{\omega}\bigl(\omega N^{1/k}+N^{2/k-1}\bigr) L^3
+N^{1/k} L^3
+
 L^3
\int_{\omega}^{1/2} 
\frac{\xi  N^{1/k}+N^{2/k-1}}{\xi^2}\dx \xi
\\&
\ll_{k}
N^{1/k} L^3 \vert \log (2\omega) \vert
+
\frac{ N^{2/k-1}}{\omega} L^3
\ll_{\ell}
\frac{ N^{2/k-1}}{\omega} L^3
\end{align*}
since $N^{-c}\le \omega \le N^{1/k-1}/L$.
A similar computation proves the result in $[-1/2,-\omega]$ too.
The estimate on $T_{k}(\alpha)$ can be obtained analogously.
\end{Proof}

 \section{Proof of Theorem \ref{thm-uncond-HL}}
 
Assume $\ell_1,\ell_2\ge 2$, $\lambda<1$, where  $\lambda$ is defined in 
 \eqref{density-def-c-def}. We'll see at the end
of the proof how the conditions in the statement of this theorem
follow; remark that in this case we cannot interchange 
the role of  $\ell_1,\ell_2$. Assume 
\begin{equation}
\label{B-def-th3}
B=B(N,\eps)=
N^{\eps},
\end{equation} 
and let $H>2B$. We have 
\begin{align}
\notag
&\sum_{n=N+1}^{N+H} 
 r^{\prime}_{\ell_1,\ell_2}(n)= 
\int_{-1/2}^{1/2} V_{\ell_1}(\alpha)  T_{\ell_2}(\alpha)  U(-\alpha,H)e(-N\alpha) \, \dx \alpha
\\
\label{dissect-HL}
&=
\int_{-B/H}^{B/H} V_{\ell_1}(\alpha) T_{\ell_2}(\alpha)  U(-\alpha,H)e(-N\alpha) \, \dx \alpha 
 + \!\!
\int_{I(B/H)} 
V_{\ell_1}(\alpha) T_{\ell_2}(\alpha)  U(-\alpha,H)e(-N\alpha) \, \dx \alpha,
\end{align}
where  $I(B/H):=[-1/2,-B/H]\cup  [B/H, 1/2]$.
By  \eqref{UH-estim}, the Cauchy-Schwarz inequality 
and Lemma \ref{tolev-coda} we have
\begin{align} 
\int_{I(B/H)} &
V_{\ell_1}(\alpha)  T_{\ell_2}(\alpha)  U(-\alpha,H)e(-N\alpha) \, \dx \alpha
\notag
\ll
\Bigl(
\int_{I(B/H)} 
\vert V_{\ell_1}(\alpha) \vert ^{2} \frac{\dx \alpha}{\vert\alpha\vert}
\Bigr)^{1/2}
\Bigl(
\int_{I(B/H)} 
\vert T_{\ell_2}(\alpha) \vert ^{2} \frac{\dx \alpha}{\vert\alpha\vert}
\Bigr)^{1/2}
\\&
  \label{V-average-coda-HL}
 \ll_{\ell_1,\ell_2} 
  \frac{H N^{\lambda -1}L^2}{B},
\end{align}
provided that $H\gg \max( N^{1-1/\ell_1}; N^{1-1/\ell_2})BL$.
By \eqref{dissect-HL}-\eqref{V-average-coda-HL}, we get
\begin{align}
 \notag
\sum_{n=N+1}^{N+H} 
&r^{\prime}_{\ell_1,\ell_2}(n)
 =  
\int_{-B/H}^{B/H} V_{\ell_1}(\alpha) T_{\ell_2}(\alpha)  U(-\alpha,H)e(-N\alpha) \, \dx \alpha 
 +
\Odipg{\ell_1,\ell_2}{ \frac{H N^{\lambda -1}L^2}{B}}.
\end{align}

Hence, recalling \eqref{main-defs}, we obtain
\begin{align} 
 \notag
\sum_{n=N+1}^{N+H} 
&r^{\prime}_{\ell_1,\ell_2}(n)  
 = \int_{-B/H}^{B/H} f_{\ell_1}(\alpha) f_{\ell_2}(\alpha)  U(-\alpha,H)e(-N\alpha) \, \dx \alpha 
 \\
 \notag
 &
\hskip0.5cm
+
\int_{-B/H}^{B/H} f_{\ell_2}(\alpha) (V_{\ell_1}(\alpha) - f_{\ell_1}(\alpha) ) U(-\alpha,H)e(-N\alpha) \, \dx \alpha 
\\
 \notag
 &
\hskip0.5cm
+
\int_{-B/H}^{B/H} f_{\ell_1}(\alpha) (T_{\ell_2}(\alpha) - f_{\ell_2}(\alpha) ) U(-\alpha,H)e(-N\alpha) \, \dx \alpha 
 \\
 \notag
 &
 \hskip0.5cm
 +
\int_{-B/H}^{B/H} (V_{\ell_1}(\alpha)-f_{\ell_1}(\alpha) ) (T_{\ell_2}(\alpha) - f_{\ell_2}(\alpha) ) U(-\alpha,H)e(-N\alpha) \, \dx \alpha
+
\Odipg{\ell_1,\ell_2}{ \frac{H N^{\lambda -1}L^2}{B}}
 \\
 \label{main-dissection-HL}
 & 
= \I_1 +\I_2 + \I_3 + \I_4 + E,
 \end{align}
 say. We now evaluate these terms.

\subsection{Computation of the main term $\I_1$}
Recalling  
$I(B/H)=[-1/2,-B/H]\cup  [B/H, 1/2]$
and Definition \eqref{density-def-c-def}, 
a direct calculation, \eqref{UH-estim}, 
the Cauchy-Schwarz inequality
and Lemma \ref{media-f-ell} give
\begin{align}
\notag
\I_1
&
=
\sum_{n=1}^H 
\int_{-1/2}^{1/2} f_{\ell_1}(\alpha) f_{\ell_2}(\alpha)   e(-(n+N)\alpha)\, \dx \alpha 
+\Odipg{\ell_1,\ell_2}{
\int_{I(B/H)} \vert f_{\ell_1}(\alpha) f_{\ell_2}(\alpha)\vert  \frac{ \dx \alpha}{ \vert \alpha \vert} 
}
\\
\notag
&
=
\frac{1}{\ell_1\ell_2}
\sum_{n=1}^H \sum_{\substack {m_{1}+m_{2} =n+N\\ N/A \le  m_{1} \le  N \\ N/A \le  m_{2}\le  N}}
m_1^{1/\ell_1-1}m_2^{1/\ell_2-1}
+\Odipg{\ell_1,\ell_2}
{ \frac{H}{B}
F_{\ell_1}(N,A)^{1/2}
F_{\ell_2}(N,A)^{1/2}
}
 \\
\label{I1-eval-step}
 & 
=
M_{\ell_1,\ell_2}(H,N)  
+
\Odipg{\ell_1,\ell_2}{
\frac{H}{B}N^{\lambda-1} A^{1-\lambda}L
},
 \end{align}
 say.
Recalling Lemma 2.8
of Vaughan \cite{Vaughan1997} we can see that
 order of magnitude of the main term $M_{\ell_1,\ell_2}(H,N)$ is 
$ c(\ell_1,\ell_2) HN^{\lambda-1}$. 
We first complete the range of summation for $m_1$ and $m_2$ to the
interval $[1, N]$.
The corresponding error term is
\begin{align*}
  &\ll_{\ell_1, \ell_2} \!\!
  \sum_{n=1}^H
    \sum_{\substack {m_{1}+m_{2} =n+N\\ 1 \le  m_{1} \le  N/A \\ 1 \le  m_{2}\le  N}}
      m_1^{1/\ell_1-1}m_2^{1/\ell_2-1}
  \ll_{\ell_1, \ell_2} \!\!
  \sum_{n=1}^H \sum_{m=1}^{N/A} m^{1/\ell_2-1} (n+N-m)^{1/\ell_1-1} \\
  &\ll_{\ell_1, \ell_2}
  H N^{1/\ell_1-1} \sum_{m=1}^{N/A} m^{1/\ell_2-1}
  \ll_{\ell_1, \ell_2} \!\!
  H N^{\lambda-1} A^{- 1/\ell_2}.
\end{align*}
We deal with the main term $M_{\ell_1,\ell_2}(H,N)$ using Lemma~2.8 of
Vaughan \cite{Vaughan1997}, which yields the $\Gamma$ factors
hidden in $c(\ell_1,\ell_2)$:
\begin{align*}
  \frac{1}{\ell_1\ell_2}
  &
  \sum_{n=1}^H \sum_{\substack {m_{1}+m_{2} =n+N\\ 1\le  m_{1} \le  N \\ 1 \le  m_{2}\le  N}}
    m_1^{1/\ell_1-1}m_2^{1/\ell_2-1}
  =
  \frac{1}{\ell_1\ell_2}
  \sum_{n=1}^H \sum_{m=1}^{N} m^{1/\ell_2-1}(n+N-m)^{1/\ell_1-1} \\
  &=
  c(\ell_1,\ell_2)
  \sum_{n=1}^H
   \Bigl[ (n+N)^{\lambda-1} +
     \Odim{(n+N)^{1/\ell_1-1} + N^{1/\ell_2-1}n^{1/\ell_1}}
  \Bigr] \\
  &= 
  c(\ell_1,\ell_2)
  \sum_{n=1}^H (n+N)^{\lambda-1} 
  +
  \Odipm{\ell_1,\ell_2}{H N^{1/\ell_1-1} +H^{1/\ell_1+1} N^{1/\ell_2-1}} \\
  &= 
  c(\ell_1,\ell_2) HN^{\lambda-1}
  +\Odipm{\ell_1,\ell_2}{ H^2N^{\lambda-2} + H   N^{1/\ell_1-1} +H^{1/\ell_1+1}   N^{1/\ell_2-1}}.
\end{align*}
Summing up,
\begin{align}
\notag
  M_{\ell_1,\ell_2}(H,N) 
  &=
  c(\ell_1,\ell_2) HN^{\lambda-1}
  \\&
\label{main-term}
  +\Odipm{\ell_1,\ell_2}{ H^2N^{\lambda-2} + H   N^{1/\ell_1-1} +H^{1/\ell_1+1}   N^{1/\ell_2-1} 
  + HN^{\lambda-1}  A^{- 1/\ell_2}}.
\end{align}

Combining \eqref{I1-eval-step}-\eqref{main-term} and using \eqref{A-def} and \eqref{B-def-th1} we get
\begin{equation}
\label{I1-eval-HL}
\I_1
=
c(\ell_1,\ell_2) HN^{\lambda-1}
+
\Odipm{\ell_1,\ell_2}{  
HN^{\lambda-1}    A(N, -C)
},
\end{equation}
for a suitable choice of   $C=C(\eps)>0$,  provided that $H\ll N^{1-\eps}$.

\subsection{Estimate of $\I_2$}
\label{estim-I2}
Using \eqref{T-ell-estim} we obtain
\begin{equation}
\label{V-approx}
\vert V_{\ell}(\alpha) - f_{\ell}(\alpha) \vert
\le
\vert V_{\ell}(\alpha) - T_{\ell}(\alpha) \vert
+
\Odipm{\ell}{(1+\vert \alpha \vert N)^{1/2}}.
\end{equation}
Hence
 \begin{align}
\notag
\I_2
& \ll_{\ell_1}
  \int_{-B/H}^{B/H} 
\vert f_{\ell_2}(\alpha) \vert \vert V_{\ell_1}(\alpha) -T_{\ell_1}(\alpha)  \vert\vert U(-\alpha,H)\vert \, \dx \alpha
\\&
\hskip1cm
\label{I2-split}
 +
   \int_{-B/H}^{B/H} 
\vert  f_{\ell_2}(\alpha) \vert  (1+\vert \alpha \vert N)^{1/2} \vert U(-\alpha,H) \vert \, \dx \alpha
= E_1+E_2,
 \end{align}
 say.
%
%
%
Letting
\begin{align}
\notag
W(N,H,B) 
&:=
\int_{-B/H}^{B/H} (1+\vert \alpha \vert N)
  \vert U(-\alpha,H) \vert ^2\ \dx \alpha
  \ll
  \frac{H^2}{N} 
+ 
H^2 N \int_{1/N}^{1/H}\!\!\! \alpha \ \dx \alpha 
+ 
N \int_{1/H}^{B/H} \frac{\dx \alpha}{\alpha} 
\\
\label{W-def}
&\ll 
NL,
\end{align}
in which we used \eqref{UH-estim}, 
using the  Cauchy-Schwarz inequality, \eqref{W-def} and Lemma \ref{media-f-ell} 
we get
\begin{align}
\notag
E_2 
&
 \ll
 F_{\ell_2}(N,A)^{1/2} 
 W(N,H,B)^{1/2} 
 \ll_{ \ell_2}
 \Bigl( \frac{N}{A}\Bigr)^{1/\ell_2-1/2} 
( \log A )^{1/2}
(NL)^{1/2}
\\ 
 \label{E2-estim}
&
 \ll_{ \ell_2}
 N^{1/\ell_2}  A^{1/2-1/\ell_2} L^{1/2} (\log A)^{1/2},
\end{align}
where $A$ is defined in \eqref{A-def}.
Using  \eqref{UH-estim}, the Cauchy-Schwarz inequality,   \eqref{W-def}
and Lemmas  \ref{App-BCP-Gallagher-2}-\ref{media-f-ell}
we obtain
\begin{align}
\notag
E_1&\ll 
 H 
 F_{\ell_2}(N,A)^{1/2} 
 E_{\ell_1}(K,B/H)^{1/2}
\ll_{\ell_1,\ell_2}
 H \Bigl( \frac{N}{A}\Bigr)^{1/\ell_2-1/2} 
( \log A )^{1/2}
 N^{1/\ell_1-1/2}  
 A(N, -c_1/2)
 \\ 
\label{E1-estim}  
&
 \ll_{\ell_1,\ell_2}
 H N^{\lambda-1}  
 A(N, -C),
 \end{align}
for a suitable choice of   $C=C(\eps)>0$, 
provided that   $N^{1-5/(6\ell_1)+\eps}B\le H \le N^{1-\eps}$.
Summarizing,  by \eqref{density-def-c-def}, \eqref{I2-split}-\eqref{E1-estim} we obtain
that there exists $C=C(\eps)>0$ such that 
 \begin{equation}
\label{I2-estim-HL}
\I_2 \ll_{\ell_1,\ell_2}
 H N^{\lambda-1} A(N, -C),
\end{equation}
provided that  $N^{1-5/(6\ell_1)+\eps} B \le H \le N^{1-\eps}$.

 \subsection{Estimate of $\I_3$}
Using \eqref{T-ell-estim}  we obtain that
 \begin{align*}
\notag
\I_3
\ll_{\ell_2} 
   \int_{-B/H}^{B/H} 
\vert  f_{\ell_1}(\alpha) \vert  (1+\vert \alpha \vert N)^{1/2} \vert U(-\alpha,H) \vert \, \dx \alpha 
 \end{align*}
and the right hand side is similar to $E_2$ of  \S\ref{estim-I2}; hence arguing as for
 \eqref{E2-estim} we have 
 \begin{equation}
\label{I3-estim-HL}
\I_3 \ll_{\ell_1,\ell_2}
  N^{1/\ell_1}  A^{1/2-1/\ell_1} L^{1/2} (\log A)^{1/2},
\end{equation}
where $A$ is defined in \eqref{A-def}.

\subsection{Estimate of $\I_4$}
By \eqref{main-dissection} and  \eqref{V-approx} we can write
 \begin{align}
 \notag
\I_4
& \ll_{\ell_1,\ell_2}
\int_{-B/H}^{B/H}
\vert  V_{\ell_1}(\alpha) -T_{\ell_1}(\alpha)\vert  (1+\vert \alpha \vert N)^{1/2} \vert U(-\alpha,H) \vert \, \dx \alpha
\\
\label{I4-split-HL}
&\hskip1cm
+
\int_{-B/H}^{B/H}
 (1+\vert \alpha \vert N)  \vert U(-\alpha,H) \vert \, \dx \alpha
= R_1+R_2,
 \end{align}
 say.
By  the Cauchy-Schwarz inequality,
Lemma \ref{App-BCP-Gallagher-2}  and arguing as in \eqref{E2-estim}
we have
\begin{align}
 \label{R1-estim}
R_1
\ll
 E_{\ell_1}(K,B/H)^{1/2}
 W(N,H,B)^{1/2} 
 \ll_{\ell_1}
N^{1/\ell_1} A(N, -C),
\end{align}
for a suitable choice of $C=C(\eps)>0$, 
provided that   $N^{1-5/(6\ell_1)+\eps} B \le H \le N^{1-\eps}$.

Moreover  by  \eqref{UH-estim}
we get
\begin{align}
  \label{R2-estim}
R_2
 \ll 
H  \int_{-1/N}^{1/N}   \, \dx \alpha
+
H N  \int_{1/N}^{1/H}   \alpha   \, \dx \alpha 
+
N \int_{1/H}^{B/H}    \, \dx \alpha  
\ll
\frac{NB}{H}.
\end{align}
Summarizing, by \eqref{density-def-c-def} and \eqref{I4-split-HL}-\eqref{R2-estim}, we obtain
 \begin{equation}
\label{I4-estim-HL}
\I_4 \ll_{\ell_1,\ell_2} 
 H N^{\lambda-1} A(N, -C),
\end{equation}
for a suitable choice of   $C=C(\eps)>0$, 
provided that $N^{1-5/(6\ell_1)+\eps} B \le H \le N^{1-\eps}$
and $H\gg N^{1-1/\ell_2+\eps}$.

\subsection{Final words} 

Summarizing, recalling that $\ell_1, \ell_2\ge 2$,   $\lambda<1$,
by  \eqref{main-dissection-HL}-\eqref{I3-estim-HL} and \eqref{I4-estim-HL},
by   optimising the choice of $B$ as in \eqref{B-def-th3},
we have that there exists $C=C(\eps)>0$ such that 
\begin{align*}
\sum_{n=N+1}^{N+H} &
r^{\prime}_{\ell_1,\ell_2}(n)
=
c(\ell_1,\ell_2)   H N^{\lambda-1} 
+
\Odipm{\ell_1,\ell_2}{H N^{\lambda-1} A(N, -C)},
\end{align*}
uniformly for     $\max(N^{1-5/(6\ell_1)}; N^{1-1/\ell_2})N^{\eps} \le H \le N^{1-\eps}$.
Theorem \ref{thm-uncond-HL} follows.

\renewcommand{\bibliofont}{\normalsize}
  
%

  \providecommand{\bysame}{\leavevmode\hbox to3em{\hrulefill}\thinspace}
\providecommand{\MR}{\relax\ifhmode\unskip\space\fi MR }
\providecommand{\MRhref}[2]{%
  \href{http://www.ams.org/mathscinet-getitem?mr=#1}{#2}
}
\providecommand{\href}[2]{#2}

\vskip0.25cm
\noindent
\begin{tabular}{l@{\hskip 20mm}l}
Alessandro Languasco               & Alessandro Zaccagnini\\
Universit\`a di Padova     & Universit\`a di Parma\\
 Dipartimento di Matematica  & Dipartimento di Scienze Matematiche, \\
 ``Tullio Levi-Civita'' &  Fisiche e Informatiche \\
Via Trieste 63                & Parco Area delle Scienze, 53/a \\
35121 Padova, Italy            & 43124 Parma, Italy\\
{\it e-mail}: alessandro.languasco@unipd.it        & {\it e-mail}:
alessandro.zaccagnini@unipr.it 
\end{tabular}

\end{document}